\documentclass[11pt]{amsart}

\usepackage{amscd}
\usepackage{amsmath}
\usepackage{graphicx}
\usepackage{amsfonts}
\usepackage{amssymb}
\begin{document}

\title[Lie Derivations and Generalized Lie Derivations]
{Additive Lie ($\xi$-Lie) Derivations and Generalized Lie
($\xi$-Lie) Derivations on Prime Algebras}

\author{Xiaofei Qi}
\address
{Department of Mathematics\\
 Shanxi University\\
  Taiyuan 030006\\
   R. R.
of China} \email{qixf1980@126.com}
\author{Jinchuan Hou}
\address{Department of
Mathematics\\
Taiyuan University of Technology\\
 Taiyuan 030024\\
  P. R. of China}
\email{jinchuanhou@yahoo.com.cn}

\thanks{This work is partially supported by National Natural Science
Foundation of China (No. 10771157) and Research Grant to Returned
Scholars of Shanxi (2007-38).}

\begin{abstract}

The additive (generalized) $\xi$-Lie derivations  on prime algebras
are characterized. It is shown, under some suitable assumption, that
an additive map $L$ is an additive (generalized) Lie derivation if
and only if it is the sum of an additive (generalized) derivation
and an additive map from the algebra into its center vanishing all
commutators; is an additive (generalized) $\xi$-Lie derivation with
$\xi\not=1$ if and only if it is an additive (generalized)
derivation satisfying $L(\xi A)=\xi L(A)$ for all $A$. These results
are then used to characterize additive (generalized) $\xi$-Lie
derivations on several operator algebras such as Banach space
standard operator algebras and von Neumman algebras.
\end{abstract}
\maketitle

\section{Introduction}

Let $\mathcal A$ be an associative ring (or an algebra over a field
$\mathbb{F}$). Then $\mathcal A$ is a Lie ring (Lie algebra) under
the Lie product $[A,B]=AB-BA.$ Recall that an additive (linear) map
$\delta$ from $\mathcal A$ into itself is called an additive
(linear) derivation if $\delta(AB)=\delta(A)B+A\delta(B)$ for all
$A$, $B\in {\mathcal A }$. More generally, an additive (linear) map
$L$ from $\mathcal A$ into itself is called an additive (linear) Lie
derivation if $L([A,B])=[L(A),B]+[A,L(B)]$ for all $A,B \in \mathcal
A $. The questions of characterizing  Lie derivations and revealing
the relationship between Lie derivations and derivations have
received many mathematicians' attention recently (for example, see
\cite{B2, Ch1, J, M, ZH1}).

Note that an important relation associated with the Lie product is
the commutativity. Two elements $A,B\in{\mathcal A}$ are commutative
if $AB=BA$, that is, their Lie product is zero. More generally, if
$\xi\in{\mathbb F}$ is a scalar and  if $AB=\xi BA$, we say that $A$
commutes with $B$ up to the factor $\xi$. The conception of
commutativity up to a factor for pairs of operators is also
important and has been studied in the context of operator algebras
and quantum groups (ref. \cite{BBP, Ka}). Motivated by this, we
introduced an binary operation $[A,B]_\xi =AB-\xi BA$, called the
$\xi$-Lie product of $A$ and $B$, and a conception of (generalized)
$\xi$-Lie derivations in \cite{QH}. Recall that an additive (linear)
map $L: {\mathcal A} \rightarrow {\mathcal A}$ is called a $\xi$-Lie
derivation if $L([A,B]_{\xi})=[L(A),B]_{\xi}+[A,L(B)]_{\xi}$ for all
$A,B \in {\mathcal A}$; an additive (linear) map $\delta: {\mathcal
A} \rightarrow {\mathcal A}$ is called an additive (linear)
generalized $\xi$-Lie derivation if there exists an additive
(linear) $\xi$-Lie derivation $L$ from $\mathcal A$ into itself such
that $\delta([A,B]_{\xi})=\delta(A)B-\xi \delta(B)A+AL(B)-\xi BL(A)$
for all $A,B \in {\mathcal A}$,  and $L$ is called the relating
$\xi$-Lie derivation of $\delta$. These conceptions unify several
important conceptions such as (generalized) derivations,
(generalized) Jordan derivations and (generalized) Lie derivations
(see \cite{HQ1, H1}). It is clear that a (generalized) $\xi$-Lie
derivation is a (generalized) derivation if $\xi =0$; is a
(generalized) Lie derivation if $\xi=1$; is a (generalized) Jordan
derivation if $\xi=-1$. Moreover, a characterization of
(generalized) $\xi$-Lie derivations  on triangular algebras for all
possible $\xi$ is given in \cite{QH}. Note that triangular algebras
are not prime.

The purpose of the present paper is to discuss the questions of
characterizing the Lie ($\xi$-Lie) derivations and generalized Lie
($\xi$-Lie) derivations, and revealing the relationship between such
additive maps to derivations (generalized derivations) on prime
algebras. As every (generalized) $\xi$-Lie derivation is a
(generalized) derivation if $\xi =0$, we need only consider the case
that $\xi\not=0$.

Let us recall some notions and notations. Throughout this paper,
$\mathcal A$ denotes a prime algebra over a field $\mathbb F$ (i.e.
$A{\mathcal A}B=0$ implies $A=0$ or $B=0$ for any $A,\ B\in{\mathcal
A}$) with the center ${\mathcal Z}({\mathcal A})$ and maximal right
ring of quotients ${\mathcal Q}={\mathcal Q}_{mr}({\mathcal A})$.
The center $\mathcal C$ of $\mathcal Q$ is a field which is called
the extended centroid of $\mathcal A$. The central closure
${\mathcal{AC}}$ of $\mathcal A$ is the $\mathcal C$-subalgebra of
${\mathcal Q}$ generated by $\mathcal A$. An element $A\in {\mathcal
A}$ is algebraic over ${\mathcal Z}({\mathcal A})$, if there exists
a polynomial $p\in {\mathcal P}({\mathcal Z}({\mathcal A}))$ such
that $p(A)=0$, that is, there exist $Z_0$, $Z_1$, $\cdots$, $Z_n\in
{\mathcal Z}({\mathcal A})$ such that $Z_n\neq 0$ and
$p(A)=Z_0+Z_1A+\cdots Z_nA^n=0$. In this case $n={\rm deg}(p)$ is
called the degree of $p$, and $\min\{{\rm deg}(p): p(A)=0\}$ is
called  the degree of algebraicity of $A$ over ${\mathcal
Z}({\mathcal A})$, denoted by deg$(A)$. If $A$ is not algebraic over
${\mathcal Z}({\mathcal A})$, then we  write deg$(A)=\infty$. The
degree of algebraicity of $\mathcal A$ is defined as
$\mbox{deg}({\mathcal A})=\sup \{\mbox{deg}(A): A\in {\mathcal A}\}$
(Ref. \cite{Br} for details).

This paper is organized as follows. Let ${\mathcal A}$ be a prime
algebra over a field $\mathbb F$. Assume that $\xi\in{\mathbb F}$ is
a nonzero scalar and $L: {\mathcal A}\rightarrow {\mathcal A}$ is an
additive map. It is known that, if deg${\mathcal A}\geq 3$, then $L$
is an additive Lie derivation if and only if it is the sum of an
additive derivation and an additive map into its centroid vanishing
each commutator \cite{B2}; when $\mathbb F$ is of characteristic not
2, then $L$ is a Jordan derivation if and only if $L$ is an additive
derivation \cite{H2}. In Section 2, we show that, when $\mathbb F$
is of characteristic not 2 and $\mathcal A$ is unital containing a
nontrivial idempotent $P$, then $L$ is a $\xi$-Lie derivation with
$\xi\not=\pm 1$ if and only if $L$ is an additive derivation
satisfying $L(\xi A)=\xi L(A)$ for all $A\in{\mathcal A}$ (Theorem
2.1).  This result then is used to give a characterization of
additive $\xi$-Lie derivations on factor von Neumman algebras
(Theorem 2.2). For Banach space standard operator algebras, a little
more can be said. Let $\mathcal A$ be a standard operator algebra in
${\mathcal B}(X)$, i.e., $\mathcal A$ contains all finite rank
operators (note that, we do not require that $\mathcal A$ contains
the unit $I$ and is closed under norm topology), where ${\mathcal
B}(X)$ is the Banach algebra of all bounded linear operators acting
on $X$. Let $L: {\mathcal A}\rightarrow{\mathcal B}(X)$ be an
additive map. We obtain that, if $\dim X\geq 3$, then $L$ is an
additive Lie derivation if and only if $L$ is the sum of an additive
derivation on ${\mathcal A}$ and an additive map from $\mathcal A$
into ${\mathbb F}I$ annihilating each commutator; $L$ is an additive
$\xi$-Lie derivation with $\xi\not=1$ if and only if $L$ is an
additive derivation satisfying $L(\xi A)=\xi L(A)$ (Theorem 2.3).

Section 3 is devoted to characterizing the generalized $\xi$-Lie
derivations. Assume that $\mathcal A$ is unital and
$\delta:{\mathcal A}\rightarrow{\mathcal A}$ is an additive map. We
show that, if deg${\mathcal A}\geq 3$, then $\delta$ is a
generalized Lie derivation if and only if $\delta$ is the sum of an
additive generalized derivation on ${\mathcal A}$ and an additive
map from ${\mathcal A}$ into its center annihilating all
commutators; if $\mathbb F$ is of characteristic not 2, then
$\delta$ is a generalized Jordan derivation if and only if $\delta$
is an additive generalized derivation; if $\mathbb F$ is of
characteristic not 2 and $\mathcal A$ contains a nontrivial
idempotent $P$, then $\delta$ is a generalized $\xi$-Lie derivation
with $\xi\not=\pm 1$ if and only if $\delta$ is an additive
generalized derivation satisfying $\delta(\xi A)=\xi \delta(A)$ for
all $A\in{\mathcal A}$ (Theorem 3.1). As an application, a
characterization of additive generalized $\xi$-Lie derivations on
factor von Neumman algebras and Banach space standard operator
algebras is obtained (Theorem 3.2 and Theorem 3.3).

\section{Additive Lie and $\xi$-Lie derivations}

In this section, we consider the question of characterizing the
additive Lie and $\xi$-Lie derivations on prime algebras. It is
obvious that if an additive map $L$ on an algebra ${\mathcal A}$  is
the sum of an additive derivation and an additive map from $\mathcal
A$ into its center vanishing the commutators, then $L$ is a Lie
derivation. Also, it is clear that, for $\xi\not=1$, every additive
derivation $L$ satisfying $L(\xi A)=\xi L(A)$ is a $\xi$-Lie
derivation. Our main purpose in this section is to show that the
inverses of these facts are true   under some weak assumptions.

The following is the main result.

\textbf{Theorem 2.1.} {\it Let $\mathcal A$  be a prime algebra over
a field ${\mathbb F}$. Assume that $\xi\in{\mathbb F}$ is a nonzero
scalar and $L: {\mathcal A}\rightarrow {\mathcal A}$ is an additive
$\xi$-Lie derivation. }

(1) {\it If $\xi=1$, that is, if $L$ is a Lie derivation, and if
deg$\mathcal A\geq 3$, then $L(A)=\tau(A)+h(A)$ for all
$A\in{\mathcal A}$, where $\tau:{\mathcal A}\rightarrow {\mathcal
{AC}}$ (the central closure of $\mathcal A$) is an additive
derivation and $h:{\mathcal A}\rightarrow {\mathcal C}$ (the
extended centroid of $\mathcal A$) is an additive map  vanishing
each commutator.}

(2) {\it If $\xi=-1$, that is, if $L$ is a Jordan derivation, and if
$\mathbb F$ is of characteristic not 2, then $L$ is an additive
derivation.}

(3) {\it If $\xi\not=\pm 1$, $\mathbb F$ is of characteristic not 2,
 $\mathcal A$ is unital and contains a nontrivial idempotent $P$,
then $L$ is an additive derivation and satisfies $L(\xi A)=\xi L(A)$
for all $A\in{\mathcal A}$.}

{\bf Proof.} By \cite{B2},  the statement (1) is true; by \cite{H2},
the statement (2) is true.

We'll prove the statement (3) by checking  several claims. In the
sequel, we always assume that $L: {\mathcal A}\rightarrow {\mathcal
A}$ is an additive $\xi$-Lie derivation with $\xi\not=\pm 1$.

Let ${\mathcal A}_{11}=P{\mathcal A}P$, ${\mathcal
A}_{12}=P{\mathcal A}(I-P)$, ${\mathcal A}_{21}=(I-P){\mathcal A}P$
and ${\mathcal A}_{22}=(I-P){\mathcal A}(I-P)$. It is clear that
${\mathcal A}={\mathcal A}_{11}\dotplus {\mathcal A}_{12}\dotplus
{\mathcal A}_{21}\dotplus{\mathcal A}_{22}$.

\textbf{Claim 1.} $L(P)=PL(P)+(I-P)L(P)P$ and
$L(I-P)=-PL(P)(I-P)-(I-P)L(P)P+(I-P)L(I-P)(I-P)$.

Since $$\begin{array}{rl}0=&L([P,I-P]_{\xi})=[L(P),I-P]_{\xi}+[P,L(I-P)]_{\xi}\\
=&L(P)(I-P)-\xi(I-P)L(P)+PL(I-P)-\xi L(I-P)P,\end{array} \eqno
(2.1)$$ multiplying by $I-P$ from both sides in Eq.(2.1), we get
$(I-P)L(P)(I-P)-\xi(I-P)L(P)(I-P)=0$, that is,
$(1-\xi)(I-P)L(P)(I-P)=0$. Note that $\xi\not=1$. It follows that
$(I-P)L(P)(I-P)=0$. Hence
$L(P)=PL(P)P+PL(P)(I-P)+(I-P)L(P)P=PL(P)+(I-P)L(P)P$.

By Eq.(2.1), we have
$$\begin{array}{rl}0=&PL(P)(I-P)-\xi(I-P)L(P)P+PL(I-P)P\\
&+PL(I-P)(I-P)-\xi PL(I-P)P-\xi (I-P)L(I-P)P\\
=&(1-\xi)PL(I-P)P+(PL(P)(I-P)+PL(I-P)(I-P))\\
&-\xi((I-P)L(P)P+(I-P)L(I-P)P).\end{array}$$ Since $\xi\not=0,1$,
 we get $PL(I-P)P=0$,
$PL(I-P)(I-P)=-PL(P)(I-P)$ and $(I-P)L(I-P)P=-(I-P)L(P)P$.
So $$\begin{array}{rl}L(I-P)=&PL(I-P)(I-P)+(I-P)L(I-P)P+(I-P)L(I-P)(I-P)\\
=&-PL(P)(I-P)-(I-P)L(P)P+(I-P)L(I-P)(I-P).\end{array}$$

\textbf{Claim 2.} $L(I)=0$.

Define a map $L^\prime:{\mathcal A}\rightarrow{\mathcal A}$ by
$$L^\prime(A)=L(A)-[A,PL(P)(I-P)-(I-P)L(P)P]\quad{\rm for\ all\ }A\in{\mathcal A}.$$
By Claim 1, it is easy to check that
 $L^\prime$ is also an additive $\xi$-Lie derivation and satisfies that
 $$L^\prime(P)=PL(P)P\in{\mathcal
 A}_{11}\quad{\rm and}\ \  L^\prime(I-P)=(I-P)L(I-P)(I-P)\in{\mathcal A}_{22}.\eqno(2.2)$$

For any $A_{12}\in {\mathcal A}_{12}$, by Eq.(2.2), we have
$A_{12}L^\prime(P)=0$ and $L^\prime(I-P)A_{12}=0$. Since
$$L^\prime(A_{12})=L^\prime([P,A_{12}]_\xi)=L^\prime(P)A_{12}-\xi A_{12}L^\prime(P) +P L^\prime(A_{12})-\xi
L^\prime(A_{12})P,\eqno(2.3)$$ multiplying by $(I-P)$ from the right
side in Eq.(2.3), we get
$$L^\prime(A_{12})(I-P)=L^\prime(P)A_{12}(I-P)+PL^\prime(A_{12})(I-P)=L^\prime(P)A_{12}+PL^\prime(A_{12})(I-P).$$
Multiplying by $P$ from the left side in the above equation, we have
$L^\prime(P)A_{12}=0$. Hence we have proved that
$$A_{12}L^\prime(P)=L^\prime(P)A_{12}=0.\eqno(2.4)$$

Similarly, by using of the relation $L^\prime(I-P)\in{\mathcal
A}_{22}$, one can show that
$$L^\prime(I-P)A_{12}=A_{12}L^\prime(I-P)=0.\eqno(2.5)$$
Combining Eq.(2.4) with (2.5), we obtain
$L^\prime(I)A_{12}=A_{12}L^\prime(I)=0.$ Since $\mathcal A$ is
prime, it follows that
$$L^\prime(I)P=0 \quad{\rm and}\ \ (I-P)L^\prime(I)=0.\eqno(2.6)$$

Now for any $A\in{\mathcal A}$, since $\xi\not=1$ and
$0=L^\prime([A,I]_\xi)-L^\prime([I,A]_\xi)=(1-\xi)[A,L^\prime(I)],$
we get $L^\prime(I)\in {\mathcal Z}({\mathcal A})$. Thus, by
Eq.(2.6), we have $L^\prime(I)=0$, and so $L(I)=0$. Complete the
proof of the claim.

{\bf Claim 3.} For any $A\in {\mathcal A}$, we have $L(\xi A)=\xi
L(A)$ and $L$ is an additive derivation.

For any $A\in {\mathcal A}$, by the definition of $L$, we have
$$L((1-\xi)A)=L([I,A]_\xi)=L(I)A-\xi AL(I)+L(A)-\xi L(A),$$ that is,
$$-L(\xi A)=L(I)A-\xi AL(I)-\xi L(A).$$
This and Claim 2 yield to $$L(\xi A)=\xi L(A).\eqno(2.7)$$

Now take any $A,B\in {\mathcal A}$. Note that
$(1-\xi)[A,B]_{-1}=[A,B]_\xi+[B,A]_\xi$ and $\xi\not=1$. Then, by
Eq.(2.7), we have
$$\begin{array}{rl}
&L((1-\xi)[A,B]_{-1})=L([A,B]_\xi)+L([A,B]_\xi)\\
=&L(A)B-\xi BL(A)+AL(B)-\xi L(B)A+L(B)A-\xi AL(B)+BL(A)-\xi L(A)B\\
=&(1-\xi)(L(A)B+AL(B)+L(B)A+BL(A),\end{array}$$that is,
$$L(AB+BA)=L(A)B+AL(B)+L(B)A+BL(A).$$  Hence $L$ is
an additive Jordan derivation from $\mathcal A$ into itself. By
statement (2), $L$ is an additive derivation, completing the proof
of the theorem. \hfill$\square$

As an application of Theorem 2.1 to the factor von Neumman algebras
case, we have

\textbf{Theorem 2.2.} {\it Let ${\mathcal M}$ be a factor von
Neumann algebra and $\xi\in{\mathbb C}$  a nonzero scalar. Assume
that  $L: {\mathcal M}\rightarrow{\mathcal M}$ is an additive
 $\xi$-Lie derivation. }

(1) {\it If $\xi=1$ and $\deg {\mathcal M}\geq 3$, then there exist
an additive derivation $\tau$ on $\mathcal M$ and an additive
functional $h: {\mathcal M}\rightarrow{\mathbb C}$ vanishing on each
commutator such that $L(A)=\tau(A)+h(A)I$ for all $A\in{\mathcal
M}$. }

(2) {\it If $\xi\not=1$, then $L$ is an additive  derivation and
satisfies that $L(\xi A)=\xi L(A)$ for all $A\in{\mathcal M}$.}

Recall that a subalgebra ${\mathcal A}\subseteq {\mathcal B}(X)$ is
called a standard operator algebra if it contains all finite rank
operators of ${\mathcal B}(X)$. Note that $\mathcal A$ may not
contain the unit operator $I$ and Theorem 2.1 can not be applied.
For the standard operator algebra ${\mathcal A}$, we have the
following result.

\textbf{Theorem 2.3.} {\it Let $X$ be a Banach space over the real
or complex field $\mathbb F$ and  ${\mathcal A}$ a standard operator
subalgebra of $\mathcal B(X)$. Assume that $\xi\in{\mathbb F}$ with
$\xi\not=0$ and $L:{\mathcal A}\rightarrow {\mathcal B}(X)$ is an
additive $\xi$-Lie derivation.}

(1) {\it If $\xi=1$, that is, if $L$ is a Lie derivation, and if
$\dim X\geq 3$,
 then $L(A)=\tau(A)+h(A)I$ for all $A\in{\mathcal A}$, where
 $\tau:{\mathcal A}\rightarrow {\mathcal B}(X)$ is an additive derivation
 and $h:{\mathcal A}\rightarrow {\mathbb F}$ is an additive map vanishing all
 commutators.}

(2) {\it If $\xi=-1$, that is, if $L$ is a Jordan derivation,  then
$L$ is an additive derivation.}

(3) {\it If $\xi\not=\pm 1$, then $L$ is an additive derivation and
satisfies  $L(\xi A)=\xi L(A)$ for all $A\in{\mathcal A}$.}

We remark that, if $X$ is infinite dimensional, then, by \cite{P},
every additive derivation $\tau$ on ${\mathcal A}$ is in fact inner,
that is, there exists an operator $T\in{\mathcal B}(X)$ such that
$\tau(A)=TA-AT$ for all $A\in{\mathcal A}$; if $X$ is finite
dimensional, then every additive derivation $\tau$ on $M_n({\mathbb
F})$ has the form $\tau(A)=TA-AT+(f(a_{ij}))_{n\times n}$ for all
$A=(a_{ij})_{n\times n}\in M_n({\mathbb F})$, where $T\in
M_n({\mathbb F})$ and $f:{\mathbb F}\rightarrow{\mathbb F}$ is an
additive derivation.

{\bf Proof of Theorem 2.3.} By Theorem 2.1(1), the statements (1)
and (2) are true.

We'll complete the proof of the statement (3) by checking several
claims. Fix a nontrivial idempotent $P\in{\mathcal A}$. In the
sequel, as a notational convenience, we denote $\mathcal
A_{11}=P{\mathcal A}P$, $\mathcal A_{12}=\{PA-PAP : A\in{\mathcal
A}\}$, $\mathcal A_{21}=\{AP-PAP : A\in{\mathcal A}\}$ and $\mathcal
A_{22}=\{ A-AP-PA+PAP: A\in{\mathcal A}\}$. Thus ${\mathcal
A}={\mathcal A}_{11}\dot{+}{\mathcal A}_{12}\dot{+}{\mathcal
A}_{21}\dot{+}{\mathcal A}_{22}$. Similarly, write ${\mathcal
B}(X)={\mathcal B}_{11}\dot{+}{\mathcal B}_{12}\dot{+}{\mathcal
B}_{21}\dot{+}{\mathcal B}_{22}$. Assume that $\xi\not=1$ and $L:
{\mathcal A}\rightarrow{\mathcal B}(X)$ is an additive $\xi$-Lie
derivation.

\textbf{Claim 1.} {\it $PL(P)P=(I-P)L(P)(I-P)=0.$}

For any $A_{22}\in{\mathcal A}_{22}$, by the definition of $L$,  we
have
$$0=L([P,A_{22}]_\xi)=L(P)A_{22}-\xi A_{22}L(P)+PL(A_{22})-\xi L(A_{22})P.\eqno(2.8)$$
Multiplying $I-P$ from the both sides of Eq.(2.8), we get
$$(I-P)L(P)A_{22}=\xi A_{22}L(P)(I-P)\quad{\rm for \ all}\
A_{22}\in{\mathcal A}_{22}.\eqno(2.9)$$ Since ${\mathcal
F}(X)\subseteq\mathcal A$ is dense in ${\mathcal B}(X)$ under the
strong operator topology, there exists a net
$\{A_\alpha\}\subset{\mathcal F}(X)$ such that SOT-$\lim_\alpha
A_\alpha= I$. Note that $A_\alpha-PA_\alpha-A_\alpha P+PA_\alpha
P\in{\mathcal A}_{22}$ and $A_\alpha-PA_\alpha-A_\alpha P+PA_\alpha
P\rightarrow I-E$ strongly. Replacing $A_{22}$ by
$A_\alpha-PA_\alpha-A_\alpha P+PA_\alpha P$ in Eq.(2.9), we get
$(I-P)L(P)(I-P)=0$ since $\xi\not=1$.

For any $A_{12}\in \mathcal A_{12}$, we have
$$L(A_{12})=L(PA_{12}-\xi A_{12}P)
=L(P)A_{12}-\xi A_{12}L(P)+PL(A_{12})-\xi L(A_{12})P.$$ Multiplying
$I-P$ from the right side of the above equation, we get
$$L(A_{12})(I-P)=L(P)A_{12}(I-P)-\xi
A_{12}L(P)(I-P)+PL(A_{12})(I-P),$$ that is,
$$\begin{array}{rl}(I-P)L(A_{12})(I-P)=&L(P)A_{12}-\xi
A_{12}(I-P)L(P)(I-P)\\
=&PL(P)A_{12}+(I-P)L(P)A_{12}.\end{array}$$ This implies that
$PL(P)A_{12}=0$. Since $\mathcal A$ is prime, it follows that
$PL(P)P=0$, completing the proof of the claim.

Now, define a map $L^\prime:{\mathcal A}\rightarrow{\mathcal B}(X)$
by
$$L^\prime(A)=L(A)-[A,PL(P)(I-P)-(I-P)L(P)P]\quad{\rm for\ all\ }A\in{\mathcal A}.$$
By Claim 1, it is easy to check that $L^\prime$ is also an additive
$\xi$-Lie derivation and satisfies that $L^\prime(P)=0$.

The following we'll prove that $L^\prime$ is an additive derivation,
and so $L$ is an additive derivation, as desired.

\textbf{Claim 2.} {\it  $L^\prime({\mathcal
A}_{ii})\subseteq{\mathcal B}_{ii}$, $i=1,2$.}

For any $A_{22}\in{\mathcal A}_{22}$, we have
$$\begin{array}{rl}0=L^\prime([P,A_{22}]_\xi)
=&L^\prime(P)A_{22}-\xi A_{22}L^\prime(P)+PL^\prime(A_{22})-\xi
L^\prime(A_{22})P\\
=&PL^\prime(A_{22})-\xi L^\prime(A_{22})P.\end{array}$$ That is,
$$PL^\prime(A_{22})P+PL^\prime(A_{22})(I-P)-\xi PL^\prime(A_{22})P-\xi(I-P)L^\prime(A_{22})P=0.$$
Note that $\xi\not=1$. It follows that
$PL^\prime(A_{22})P=PL^\prime(A_{22})(I-P)=(I-P)L^\prime(A_{22})P=0$,
and so $L^\prime(A_{22})\in{\mathcal B}_{22}$.

Taking any $A_{11}\in{\mathcal A}_{11}$ and $A_{22}\in{\mathcal
A}_{22}$, we have
$$\begin{array}{rl}0=L^\prime([A_{11},A_{22}]_\xi)
=&L^\prime(A_{11})A_{22}-\xi
A_{22}L^\prime(A_{11})+A_{11}L^\prime(A_{22})-\xi
L^\prime(A_{22})A_{11}\\
=&L^\prime(A_{11})A_{22}-\xi A_{22}L^\prime(A_{11})\\
=&PL^\prime(A_{11})(I-P)A_{22}+(I-P)L^\prime(A_{11})(I-P)A_{22}\\
&-\xi A_{22}(I-P)L^\prime(A_{11})P-\xi
A_{22}(I-P)L^\prime(A_{11})(I-P).\end{array}$$ This implies that
$$PL^\prime(A_{11})(I-P)A_{22}=0,\ \ \xi A_{22}(I-P)L^\prime(A_{11})P=0\eqno(2.10)$$
and $$(I-P)L^\prime(A_{11})(I-P)A_{22}=\xi
A_{22}(I-P)L^\prime(A_{11})(I-P).\eqno(2.11)$$ Since ${\mathcal
F}(X)\subseteq\mathcal A$ is dense in ${\mathcal B}(X)$ under the
strong operator topology, there exists a net
$\{A_\alpha\}\subset{\mathcal F}(X)$ such that SOT-$\lim_\alpha
A_\alpha= I$. Note that $A_\alpha-PA_\alpha-A_\alpha P+PA_\alpha
P\in{\mathcal A}_{22}$ and $A_\alpha-PA_\alpha-A_\alpha P+PA_\alpha
P\rightarrow I-E$ strongly. Replacing $A_{22}$ by
$A_\alpha-PA_\alpha-A_\alpha P+PA_\alpha P$ in Eqs.(2.10)-(2.11), we
get
$PL^\prime(A_{11})(I-P)=(I-P)L^\prime(A_{11})P=(I-P)L^\prime(A_{11})(I-P)=0$
since $\xi\not=0,1$. Hence $L^\prime(A_{11})\in{\mathcal B}_{11}$,
completing the proof of the claim.

\textbf{Claim 3.} {\it  $L^\prime({\mathcal A}_{ij})\subseteq
{\mathcal B_{ij}}$, $1\leq i\not=j\leq 2$.}

For any $A_{12}\in \mathcal A_{12}$, noting that $L^\prime(P)=0$, we
have
$$\begin{array}{rl}L^\prime(A_{12})=&L^\prime(PA_{12}-\xi A_{12}P)\\
=&L^\prime(P)A_{12}-\xi A_{12}L^\prime(P)+PL^\prime(A_{12})-\xi
L^\prime(A_{12})P\\
=&PL^\prime(A_{12})-\xi L^\prime(A_{12})P.\end{array}\eqno(2.12)$$
Multiplying $P$ from both sides of the above equation, we get $\xi
PL^\prime(A_{12})P=0$, which implies that $PL^\prime(A_{12})P=0$.
Similarly, multiplying $I-P$ from the left side of Eq.(2.12) leads
to
$$(I-P)L^\prime(A_{12})=-\xi (I-P)L^\prime(A_{12})P.\eqno(2.13)$$
Multiplying $I-P$ from the right side of Eq.(2.13), we get
$(I-P)L^\prime(A_{12})(I-P)=0$. Multiplying $P$ from the right side
of Eq.(2.13), we get $(1+\xi)(I-P)L^\prime(A_{12})P=0$, which
implies that $(I-P)L^\prime(A_{12})P=0$ since $\xi\not=-1$.

Similarly, for any $A_{21}\in \mathcal A_{21}$, by using of the
equation $L^\prime(A_{21})=L^\prime([A_{21},P]_\xi)$, one can check
that $PL^\prime(A_{21})P=0$, $(I-P)L^\prime(A_{21})(I-P)=0$ and
$PL^\prime(A_{21})(I-P)=0$.

Thus we obtain $L^\prime(A_{ij})\in {\mathcal B_{ij}}$ with
$i\not=j$.

{\bf Claim 4.} {\it  $L^\prime$ has the following properties:}

(a) {\it
$L^\prime(A_{ii}B_{ij})=L^\prime(A_{ii})B_{ij}+A_{ii}L^\prime(B_{ij})$
holds for all $ A_{ii}\in \mathcal A_{ii}$ and $ B_{ij}\in \mathcal
A_{ij}$, $1\leq i\not=j\leq 2$.}

(b) {\it
$L^\prime(A_{ij}B_{jj})=L^\prime(A_{ij})B_{jj}+A_{ij}L^\prime(B_{jj})$
holds for all $A_{ij}\in{\mathcal A}_{ij}$ and $ B_{jj}\in{\mathcal
A}_{jj}$, $1\leq i\not=j\leq 2$.}

(c) {\it
$L^\prime(A_{ij}B_{ji})=L^\prime(A_{ij})B_{ji}+A_{ij}L^\prime(B_{ji})$
 holds for all  $A_{ij}\in \mathcal A_{ij}$ and $ B_{ji}\in \mathcal A_{ji}$,
$1\leq i\not=j\leq 2$.}

(d) {\it
$L^\prime(A_{ii}B_{ii})=L^\prime(A_{ii})B_{ii}+A_{ii}L^\prime(B_{ii})$
 holds for all  $A_{ii},B_{ii} \in \mathcal A_{ii}$, $i=1,2$.}

For any $A_{ii}\in \mathcal A_{ii}$ and $B_{ij}\in \mathcal A_{ij},$
it follows from Claims 2-3 that
$$\begin{array}{rl}L^\prime(A_{ii}B_{ij})
=&L^\prime([A_{ii},B_{ij}]_\xi)\\
=&L^\prime(A_{ii})B_{ij}-\xi
B_{ij}L^\prime(A_{ii})+A_{ii}L^\prime(B_{ij})-\xi
L^\prime(B_{ij})A_{ii}\\
=&L^\prime(A_{ii})B_{ij}+A_{ii}L^\prime(B_{ij}),\end{array}$$ and so
(a) holds true.

Similarly, (b) is true for all $A_{ij}\in \mathcal A_{ij}$ and
$B_{jj}\in \mathcal A_{jj}.$

For any $A_{ij}\in \mathcal A_{ij}$ and $ B_{ji}\in \mathcal
A_{ji}$, by Claim 3 and the additivity of $L^\prime$, we get
$$\begin{array}{rl}L^\prime(A_{ij}B_{ji})-L^\prime(\xi B_{ji}A_{ij})
=&L^\prime([A_{ij},B_{ji}]_\xi)\\
=&L^\prime(A_{ij})B_{ji}-\xi
B_{ji}L^\prime(A_{ij})+A_{ij}L^\prime(B_{ji})-\xi
L^\prime(B_{ji})A_{ij}\\
=&(L^\prime(A_{ij})B_{ji}+A_{ij}L^\prime(B_{ji}))-\xi(B_{ji}L^\prime(A_{ij})+L^\prime(B_{ji})A_{ij}).\end{array}$$
Note that $L^\prime(A_{ij}B_{ji}),
L^\prime(A_{ij})B_{ji}+A_{ij}L^\prime(B_{ji})\in{\mathcal B}_{ii}$
and $L^\prime(\xi B_{ji}A_{ij}),
B_{ji}L^\prime(A_{ij})+L^\prime(B_{ji})A_{ij}\in{\mathcal B}_{jj}$.
It follows that
$L^\prime(A_{ij}B_{ji})=L^\prime(A_{ij})B_{ji}+A_{ij}L^\prime(B_{ji})$,
and so (c) is true.

For any $A_{ii},B_{ii}\in \mathcal A_{ii}$ and any
$C_{ij}\in{\mathcal A}_{ij}$, by (a), we have
$$\begin{array}{rl}L^\prime(A_{ii}B_{ii}C_{ij})=&L^\prime(A_{ii})B_{ii}C_{ij}+A_{ii}L^\prime(B_{ii}C_{ij})\\
=&L^\prime(A_{ii})B_{ii}C_{ij}+A_{ii}L^\prime(B_{ii})C_{ij}+A_{ii}B_{ii}L^\prime(C_{ij})\end{array}$$
and
$$L^\prime(A_{ii}B_{ii}C_{ij})=L^\prime(A_{ii}B_{ii})C_{ij}+A_{ii}B_{ii}L^\prime(C_{ij}).$$
Comparing the above two equations gives
$$(L^\prime(A_{ii}B_{ii})-L^\prime(A_{ii})B_{ii}-A_{ii}L^\prime(B_{ii}))C_{ij}=0$$
for all $C_{ij}\in{\mathcal A}_{ij}$. Since $\mathcal A$ is prime,
it follows that
$L^\prime(A_{ii}B_{ii})-L^\prime(A_{ii})B_{ii}-A_{ii}L^\prime(B_{ii})=0$,
that is, (d) holds true.

{\bf Claim 5.} {\it $L^\prime$ is an additive derivation, and
therefore, $L$ is an additive derivation and satisfies $L(\xi A)=\xi
L(A)$ for all $A\in{\mathcal A}$.}

For any $A,B\in{\mathcal A}$, write $A=A_{11}+A_{12}+A_{21}+A_{22}$
and $B=B_{11}+B_{12}+B_{21}+B_{22}$. By Claim 4 and the additivity
of $L^\prime$, it is easily checked that
$L^\prime(AB)=L^\prime(A)B+AL^\prime(B)$, that is, $L^\prime$ is an
additive derivation on ${\mathcal A}$. Note that the map $A\mapsto
[A,PL(P)(I-P)-(I-P)L(P)P]$ is an inner derivation of $\mathcal A$.
So $L$ is also an additive derivation.

Finally, for any $A,B\in{\mathcal A}$, we have
$$\begin{array}{rl} [L(A),B]_\xi +[A,L(B)]_\xi =& L([A,B]_\xi)=L(AB)-L(\xi
BA)\\ =& L(A)B+AL(B)-L(B)(\xi A)-BL(\xi A),\end{array}$$ which
implies that $$BL(\xi A)=\xi BL(A) \eqno(2.14)$$ holds for all
$A,B\in{\mathcal A}$. Taking a net $\{B_\alpha\}$ in $\mathcal A$
such that $B_\alpha\rightarrow I$ strongly, and replacing $B$ by
$B_\alpha$ in Eq.(2.14), we obtain $L(\xi A)=\xi L(A)$. This
completes the proof of the statement (3) in Theorem 2.3.
\hfill$\Box$

\section{Additive generalized Lie and  $\xi$-Lie derivations}

In this section, we discuss the question of characterizing the
additive generalized Lie derivations and generalized $\xi$-Lie
derivations. It is obvious that, for $\xi\not=1$, every additive
generalized derivation $\delta$ satisfying $\delta(\xi A)=\xi\delta
(A)$ is an additive generalized $\xi$-Lie derivation; and the sum of
an additive generalized derivation and an additive map into the
center vanishing all commutators is an additive generalized Lie
derivation. We show that the inverses of above facts are true for
most prime algebras.

The following is the main result in this section.

{\bf Theorem 3.1.} {\it  Let $\mathcal A$  be a unital prime algebra
over a field ${\mathbb F}$ and $\xi\in\mathbb F$ with  $\xi\not=0$.
Suppose that $\delta:{\mathcal A}\rightarrow{\mathcal A}$ is an
additive generalized $\xi$-Lie derivation with
 $L:{\mathcal A}\rightarrow{\mathcal A}$ the relating $\xi$-Lie derivation.
}

(1) {\it If $\xi=1$, that is, if $\delta$ is a generalized Lie
derivation, and if deg$\mathcal A\geq 3$, then
$\delta(A)=\delta^\prime(A)+h(A)$ for all $A\in{\mathcal A}$, where
$\delta^\prime:{\mathcal A}\rightarrow {\mathcal {AC}}$  is an
additive generalized derivation and $h:{\mathcal A}\rightarrow
{\mathcal C}$ is an additive map vanishing each commutator.}

(2) {\it If $\xi=-1$, that is, if $\delta$ is a generalized Jordan
derivation, and if $\mathbb F$ is of characteristic not 2,  then
$\delta$ is an additive generalized derivation.}

(3) {\it If $\xi\not=\pm 1$, $\mathbb F$ is of characteristic not 2,
and if $\mathcal A$ contains a nontrivial idempotent, then $\delta$
is an additive generalized derivation and $\delta(\xi A)=\xi
\delta(A)$ for all $A\in{\mathcal A}$.}

{\bf Proof.} Since $\delta:{\mathcal A}\rightarrow{\mathcal A}$ is
an additive generalized $\xi$-Lie derivation with
 $L:{\mathcal A}\rightarrow{\mathcal A}$ the relating $\xi$-Lie derivation,
 we have $$\delta([A,B]_{\xi})=\delta(A)B-\xi
\delta(B)A+AL(B)-\xi BL(A)$$ for all $A,B \in {\mathcal A}$. Taking
$B=I$ in the above equation, we get $\delta(A-\xi
A)=\delta(A)-\xi\delta(I)A+AL(I)-\xi L(A)$, that is,
$$\delta(-\xi A)=-\xi L(A)-\xi\delta(I)A+AL(I) {\rm \ \ for \ all }\ \ A\in{\mathcal A}.\eqno(3.1)$$

If $\xi=1$, then Eq.(3.1) becomes $\delta(A)= L(A)+\delta(I)A+AL(I)$
for all $A\in{\mathcal A}$. By \cite{B2}, $L$ has the form of
$L(A)=\tau(A)+h(A)$, where $\tau$ is an additive derivation of
${\mathcal A}$ and $h:{\mathcal A}\rightarrow{\mathcal C}$ is an
additive map satisfying $h([A,B])=0$ for all $A$ and $B$. Define
$\delta^{\prime}:{\mathcal A}\rightarrow{\mathcal A}$ by
$\delta^{\prime}(A)=\tau(A)+\delta(I)A+AL(I)$ for all $A\in{\mathcal
A}$. Thus we get $\delta(A)=\delta^{\prime}(A)+h(A)$. It is easily
seen that $\delta^{\prime}$ is an additive generalized derivation.
Hence the statement (1) of Theorem 3.1 holds true.

If $\xi\not=1$, then, substituting $A$ by $-\xi^{-1}A$ in Eq.(3.1),
one gets $$\delta(A)=-\xi L(-\xi^{-1}A)+\delta(I)A-\xi^{-1}AL(I)
\eqno(3.2)$$ for all $A\in{\mathcal A}$. Since $L$ is an additive
$\xi$-Lie derivation, by Theorem 2.1(2) and (3), we see that $L$ is
an additive derivation satisfying $L(\xi A)=\xi L(A)$ for all $A$.
It follows from Eq.(3.2) that $\delta
(A)=L(A)+\delta(I)A-\xi^{-1}AL(I)$, which is a generalized
derivation. Furthermore, $\delta(\xi A)=L(\xi A)+\delta(I)\xi
A-\xi^{-1}\xi AL(I)=\xi L(A)+\delta(I)\xi A-\xi^{-1}\xi
AL(I)=\xi\delta(A)$. Hence, the statement (2) of Theorem 3.1 is
true. \hfill$\Box$

For the von Neumman algebra case, we have

{\bf Theorem 3.2.}  {\it Let ${\mathcal M}$ be a factor von Neumann
algebra and $\xi\in{\mathbb C}$  a nonzero scalar. Assume that
$\delta: {\mathcal M}\rightarrow{\mathcal M}$ be an additive
generalized $\xi$-Lie derivation. }

(1) {\it If $\xi=1$ and $\deg {\mathcal M}\geq 3$, then there exist
an additive generalized derivation $\tau$ on $\mathcal M$ and an
additive functional $h: {\mathcal M}\rightarrow{\mathbb C}$
vanishing on each commutator such that $\delta(A)=\tau(A)+h(A)I$ for
all $A\in{\mathcal M}$. }

(2) {\it If $\xi\not=1$, then $\delta$ is an additive generalized
derivation and $\delta(\xi A)=\xi \delta(A)$ for all $A\in{\mathcal
M}$.}

For Banach space standard operator algebras, we have

{\bf Theorem 3.3.} {\it Let $X$ be a Banach space over the real or
complex field $\mathbb F$ and  ${\mathcal A}$ a standard operator
subalgebra of $\mathcal B(X)$ containing the identity $I$. Assume
that $\xi\in{\mathbb F}$ with $\xi\not=0$ and $\delta:{\mathcal
A}\rightarrow {\mathcal B}(X)$ is an additive generalized $\xi$-Lie
derivation.}

(1) {\it If $\xi=1$ and $\dim X\geq 3$,
 then $\delta(A)=\tau(A)+h(A)I$ for all $A\in{\mathcal A}$, where
 $\tau:{\mathcal A}\rightarrow {\mathcal B}(X)$ is an additive generalized derivation
 and $h:{\mathcal A}\rightarrow {\mathbb F}$ is an additive map vanishing all
 commutators.}

(2) {\it  If $\xi\not=1$, then $\delta$ is an additive generalized
derivation and satisfies  $\delta(\xi A)=\xi \delta(A)$ for all
$A\in{\mathcal A}$.}


\end{document}